\input amstex
\input amsppt.sty

\magnification1200

\hsize13cm
\vsize19cm

\TagsOnRight

\def\StemAE{7}
\def\SlatAC{6}
\def\MeWaAA{5}
\def\KratBN{4}
\def\BailAA{3}
\def\BacRAA{2}
\def\BacRAZ{1}

\def\si{\sigma}
\def\om{\omega}
\def\3{\ss}
\def\({\left(}
\def\){\right)}
\def\coef#1{\langle#1\rangle}
\def\po#1#2{(#1)_#2}

\topmatter
\title Evaluations of some determinants of matrices
related to the Pascal triangle
\endtitle
\author C.~Krattenthaler$^\dagger$
\endauthor
\affil
Institut f\"ur Mathematik der Universit\"at Wien,\\
Strudlhofgasse 4, A-1090 Wien, Austria.\\
e-mail: KRATT\@Ap.Univie.Ac.At\\
WWW: \tt http://www.mat.univie.ac.at/\~{}kratt
\endaffil
\address Institut f\"ur Mathematik der Universit\"at Wien,
Strudlhofgasse 4, A-1090 Wien, Austria.
\endaddress
\thanks{$^\dagger$ Research partially supported by the Austrian
Science Foundation FWF, grants P12094-MAT and P13190-MAT,
and by EC's IHRP Programme,
grant RTN2-2001-00059.}\endthanks
\abstract
We prove several evaluations of determinants
of matrices, the entries of which are given by the recurrence
$a_{i,j}=a_{i-1,j}+a_{i,j-1}$, or variations thereof. These
evaluations were either conjectured or extend conjectures
by Roland Bacher [J. Th\'eorie Nombres Bordeaux {\bf 13} (2001), to appear].
\endabstract
\endtopmatter
\document

\leftheadtext{C. Krattenthaler}
\rightheadtext{Evaluations of some determinants}

\subhead 1. Introduction\endsubhead
In the recent preprint \cite{\BacRAZ} (see \cite{\BacRAA} for the
final version), Bacher considers determinants of matrices, the
entries of which are given by the Pascal triangle recurrence
$a_{i,j}=a_{i-1,j}+a_{i,j-1}$ or by similar
recurrences, or are related in some other way to such bivariate
recurrent sequences. The preprint \cite{\BacRAZ} contains results and
conjectures on, on the one hand,
 evaluations of some special determinants of that
type, and, on the other hand, recurrences satisfied by such determinants.

The purpose of this paper is to provide proofs for all the
conjectured determinant evaluations in \cite{\BacRAZ}, and, in some
cases, to prove in fact generalizations of conjectures in \cite{\BacRAZ}.
(Some of these results are stated without proof in \cite{\BacRAA}.)

In Section~2 we provide the proofs for the three conjectures in
\cite{\BacRAZ} on determinants of matrices, the entries of which are
given by the recurrence
$a_{i,j}=a_{i-1,j}+a_{i,j-1}+x\,a_{i-1,j-1}$. Section~3 contains
three theorems on binomial determinants which generalize (conjectural)
determinant evaluations
in Section~5 of \cite{\BacRAZ}. The methods that we use
to prove our theorems are, on the one hand, the ``LU-factorization
method" (see \cite{\KratBN, Sec.~2.6}) and, on the other hand, the
``identification of factors method" (see \cite{\KratBN, Sec.~2.4}),
except for one case where elementary row and column operations
suffice.

\subhead 2. Some determinants of matrices with entries given by an
extension of the Pascal triangle recurrence\endsubhead
Our first theorem proves Conjecture~1.8 in
\cite{\BacRAZ}. It is stated (without proof) as Theorem~1.5 in
\cite{\BacRAA}.

\proclaim{Theorem 1}Let $(a_{i,j})_{i,j\ge0}$ be the sequence given by the
recurrence
$$a_{i,j}=a_{i-1,j}+a_{i,j-1}+x\,a_{i-1,j-1},\quad \quad i,j\ge1
\tag2.1$$
and the initial conditions $a_{i,0}=\rho^i$ and $a_{0,i}=\si^i$,
$i\ge0$. Then
$$\det_{0\le i,j\le n-1}(a_{i,j})=(1+x)^{\binom
{n-1}2}(x+\rho+\si-\rho\si)^{n-1}.$$
\endproclaim
\demo{Proof}Let $F(u,v)=\sum _{i,j\ge0} ^{}a_{i,j}u^iv^j$ be the
bivariate generating function of the sequence $(a_{i,j})$. By
multiplying both sides of (2.1) by $u^iv^j$ and summing over all $i$ and
$j$, we get an equation for $F(u,v)$ with solution
$$F(u,v)=\frac {\frac {1-u} {1-\rho u}+\frac {1-v} {1-\si v}-1}
{1-u-v-uvx}.\tag2.2$$

Now we use the {\it LU-factorization method\/} (see
\cite{\KratBN, Sec.~2.6}). Let $M$ denote the matrix
$(a_{i,j})_{0\le i,j\le n-1}$. I claim that
$$M\cdot U=L,$$
where $U=(U_{i,j})_{0\le i,j\le n-1}$ with
$$U_{i,j}=\cases 1&i=j\\
(-1)^{j-i}\(\binom {j-1}i\si+\binom {j-1}{i-1}\)&i<j,\\
0&i>j,\endcases$$
and $L=(L_{i,j})_{0\le i,j\le n-1}$ with
$$L_{i,j}=\cases \rho^i&j=0,\\
(x+\rho+\si-\rho\si)(1+x)^{j-1}\sum _{\ell=0} ^{i-1}\binom
{\ell}{j-1}\rho^{i-\ell-1}&j>0.\endcases$$
The matrix $U$ is an upper triangular matrix with 1's on the
diagonal, whereas $L$ is a lower triangular matrix with diagonal
entries
$$1,\ (x+\rho+\si-\rho\si),\ (x+\rho+\si-\rho\si)(1+x),\
\dots,\ (x+\rho+\si-\rho\si)(1+x)^{n-2}.$$
It is obvious that the claimed factorization of $M$ immediately
implies the validity of the theorem.

For the proof of the claim we compute the
$(i,j)$-entry of
$M\cdot U$. By definition, it is
$\sum _{k=0} ^{j}a_{i,k}U_{k,j}$. Let first $j=0$.
Then the sum is equal to $a_{i,0}=\rho^i$, in accordance with the
definition of $L_{i,0}$. If $j>0$, then (here, and in the following,
$\coef{z^N}f(z)$ denotes the coefficient of $z^N$ in
the formal power series $f(z)$)
$$\align (M\cdot U)_{i,j}&=\sum _{k=0} ^{j}a_{i,k}U_{k,j}\\
&=M_{i,j}+
\sum _{k=0} ^{j-1}a_{i,k}U_{k,j}\\
&=\coef{u^iv^j}F(u,v)\(1+\sum _{k=0} ^{j-1}(-1)^{j-k}v^{j-k}
\(\binom {j-1}k\si+\binom {j-1}{k-1}\)\)\\
&=\coef{u^iv^j}F(u,v)\((1-v)^{j-1}-\si v(1-v)^{j-1}\)\\
&=\coef{u^iv^j}F(u,v)(1-v)^{j-1}(1-\si v).
\tag2.3
\endalign$$
Now we rewrite $F(u,v)$, as given by (2.2), slightly,
$$\align F(u,v)&=\frac {1-u-v+uv(x+\rho+\si-\rho\si)-xuv}
{(1-\rho u)(1-\si v)(1-u-v-uvx)}\\
&=\frac 1
{(1-\rho u)(1-\si v)}+
\frac {uv(x+\rho+\si-\rho\si)}
{(1-\rho u)(1-\si v)(1-u-v-uvx)}\\
&=\frac 1
{(1-\rho u)(1-\si v)}+
\frac {uv(x+\rho+\si-\rho\si)}
{(1-\rho u)(1-\si v)(1-u)(1-v)\(1-\frac {uv(x+1)}{(1-u)(1-v)}\)}.
\endalign$$
This is substituted in (2.3):
$$\multline (M\cdot U)_{i,j}\\=
\coef{u^iv^j}(1-v)^{j-1}
\(\frac 1
{(1-\rho u)}+
\frac {uv(x+\rho+\si-\rho\si)}
{(1-\rho u)(1-u)(1-v)\(1-\frac {uv(x+1)}{(1-u)(1-v)}\)}\).
\endmultline$$
Since the degree of $(1-v)^{j-1}$ in $v$ is smaller than $j$, it is
obvious that the first expression in parentheses does not contribute
anything to the coefficient of $u^iv^j$. Hence we obtain
$$\align (M\cdot U)_{i,j}&=
\coef{u^iv^j}(1-v)^{j-1}
\frac {uv(x+\rho+\si-\rho\si)}
{(1-u)(1-v)}\sum _{h=0} ^{\infty}\frac {u^hv^h(x+1)^h}
{(1-u)^h(1-v)^h}\sum _{\ell=0} ^{\infty}\rho^\ell u^\ell\\
&=(x+\rho+\si-\rho\si)\sum _{h,\ell\ge0} ^{}\rho^\ell(1+x)^h
\binom {i-\ell-1}h\binom {0}{h+1-j}.
\endalign$$
Because of the second binomial coefficient, the summand is
nonvanishing only for $h=j-1$. This yields exactly the claimed
expression for $L_{i,j}$.\quad \quad \qed
\enddemo

The next theorem proves the conjecture about the evaluation of
$\det(A(2n))$ on page~4 of
\cite{\BacRAZ}. It is stated (without proof) in
\cite{\BacRAA, equation for $\det(B(2n))$ after
Theorem~1.5}.

\proclaim{Theorem 2}Let $(a_{i,j})_{i,j\ge0}$ be the sequence given by the
recurrence
$$a_{i,j}=a_{i-1,j}+a_{i,j-1}+x\,a_{i-1,j-1},\quad \quad i,j\ge1
$$
and the initial conditions $a_{i,i}=0$, $i\ge0$, $a_{i,0}=\rho^{i-1}$ and
$a_{0,i}=-\rho^{i-1}$, $i\ge1$. Then
$$\det_{0\le i,j\le 2n-1}(a_{i,j})=(1+x)^{2(n-1)^2}
(x+\rho)^{2n-2}.$$
\endproclaim

\demo{Proof}
We compute again the generating function
$F(u,v)=\sum _{i,j\ge0} ^{}a_{i,j}u^iv^j$. This time we have
$$F(u,v)=\frac {(1-u-v+\rho uv)} {(1-u-v-xuv)}
\frac {(u-v)} {(1-\rho u)(1-\rho v)}.$$

Also here, we apply the LU-factorization method. However, we cannot
apply LU-factorization directly to the original matrix since it is
skew-symmetric and, hence, all
odd-dimensional principal minors vanish. Instead,
we apply LU-factorization to
the matrix that results from the original matrix by first reversing the
order of rows and then transposing the resulting matrix. This is the matrix
$M=(a_{2n-1-j,i})_{0\le i,j\le 2n-1}$.
I claim that
$$M\cdot U=L,$$
where $U=(U_{i,j})_{0\le i,j\le n-1}$ with
$$U_{i,j}=\cases 0&\kern-4.1cm i>j,\\
\frac {(-1)^{j-i}} {\rho}\(\binom {j-1}i+\binom {j-1}{i-1}\rho\)&
\kern-4.1cm i\le j<2n-1,\\
1&\kern-4.1cm i=j=2n-1,\\
-\frac {1} {(\rho+x)(1+x)^{n-2}}\Big(\sum _{s=0} ^{n-i/2-2}2^{2s+1}
\binom {n-\frac {i} {2}-1+s}{2s+1} \sum _{t=0} ^{i/2-s-1}\binom
{n-1}tx^{\frac {i} {2}-s-t-1}\\
\kern1cm
+\rho \sum _{s=0} ^{n-i/2-1}2^{2s}
\binom {n-\frac {i} {2}-1+s}{2s} \sum _{t=0} ^{i/2-s-1}\binom
{n-1}tx^{\frac {i} {2}-s-t-1}\Big)\\
&\kern-4.1cm
i\text { even, }i<j=2n-1,\\
\frac {1} {(\rho+x)(1+x)^{n-2}}\Big(\sum _{s=0} ^{n-i/2-3/2}2^{2s}
\binom {n-\frac {i} {2}-\frac {3} {2}+s}{2s} \sum _{t=0} ^{i/2-s-1/2}\binom
{n-1}tx^{\frac {i} {2}-s-t-\frac {1} {2}}\\
\kern1cm
+\rho \sum _{s=0} ^{n-i/2-3/2}2^{2s+1}
\binom {n-\frac {i} {2}-\frac {1} {2}+s}{2s+1} \sum _{t=0}
^{i/2-s-3/2}\binom
{n-1}tx^{\frac {i} {2}-s-t-\frac {3} {2}}\Big)\\
&\kern-4.1cm
i\text { odd, }i<j=2n-1,
\endcases$$
and where $L$ is a lower triangular matrix with diagonal entries
$$\multline
\rho^{2n-2},\ -\frac {\rho+x} {\rho},\ \frac {(\rho+x)(1+x)} {\rho},\
-\frac {(\rho+x)(1+x)^2} {\rho},\ \dots,\\
\frac {(\rho+x)(1+x)^{2n-3}} {\rho},\
-(1+x)^{n-1}
\endmultline\tag2.4$$
Again, it is immediately obvious that the claimed factorization of $M$
implies the validity of the theorem.

The proof of the claim requires again some calculations, which turn
out to be more tedious here.
The $(i,j)$-entry of
$M\cdot U$ is $\sum _{k=0} ^{j}a_{2n-1-k,i}U_{k,j}$. It suffices to
consider the case where $0\le i\le j\le 2n-1$.

\smallskip
{\it Case 1: $i=j=0$}. Then
$\sum _{k=0} ^{j}a_{2n-1-k,i}U_{k,j}=a_{2n-1,0}=\rho^{2n-2}$, in
accordance with (2.4).

\smallskip
{\it Case 2: $0=i<j<2n-1$}. We calculate the sum in question:
$$\align (M\cdot U)_{i,j}&=\sum _{k=0} ^{j}a_{2n-1-k,i}U_{k,j}\\
&=\coef{u^{2n-1}v^i}\sum _{k=0} ^{j}\frac {(-1)^{j-k}} {\rho}
\(\binom {j-1}k+\binom {j-1}{k-1}\rho\)u^kF(u,v)\\
&=\coef{u^{2n-1}v^i}
\frac {(-1)^{j}} {\rho}\((1-u)^{j-1}-\rho u(1-u)^{j-1}\)F(u,v)\\
&=-\coef{u^{2n-1}v^i}\frac {(u-1)^{j-1}} {\rho}
\frac {(1-u-v+\rho uv)(u-v)} {(1-u-v-xuv)(1-\rho v)}.
\tag2.5\endalign$$
If $i=0$ then the expression in the last row reduces to
$$-\coef{u^{2n-1}}\frac {(u-1)^{j-1}} {\rho}\cdot u.
$$
The degree of the polynomial in $u$ of which the coefficient of $u^{2n-1}$
has to be read off is $j$. Since $j<2n-1$, this coefficient is 0.

\smallskip
{\it Case 3: $0<i\le j<2n-1$}. We start with the expression (2.5) and
apply the partial fraction expansion (with respect to $v$)
$$\multline
\frac {(1-u-v+\rho uv)(u-v)} {(1-u-v-xuv)(1-\rho v)}
=\frac {1-\rho u} {\rho(1+xu)}-\frac {(1-\rho)(1-\rho u)}
{\rho(1-\rho v)(1-\rho+\rho u+xu)}\\
+\frac {u(\rho+x)(u^2x+2u-1)} {(1+xu)(1-\rho+\rho u+xu)\(1-\frac
{v(1+xu)} {1-u}\)}.
\endmultline\tag2.6$$
After having substituted this in (2.5), we can read off the
coefficient of $v^i$. In particular, the first term on the right-hand
side of (2.6) does not contribute anything, since we have $i>0$ by assumption.
We obtain
$$\align (M\cdot U)_{i,j}&=
-\coef{u^{2n-1}}\frac {(u-1)^{j-1}} {\rho(1-\rho+\rho u+xu)}
\bigg(-(1-\rho)(1-\rho u)
\rho^{i-1}\\
&\kern1cm
+{u(\rho+x)(u^2x+2u-1)}
\frac {(1+xu)^{i-1}} {(1-u)^i}\bigg).
\tag2.7\endalign$$

Let first $i<j$. The term $(1-u)^i$ in the denominator cancels with
the term $(u-1)^{j-1}$. But the term $(1-\rho+\rho
u+xu)$ in the denominator cancels as well. because the expression in
parentheses contains this term as a factor. This can be seen by
observing that the expression in parentheses vanishes for
$u=(\rho-1)/(\rho+x)$ (which is the zero of the term
$(1-\rho+\rho u+xu)$). Hence, the task in (2.7) is to read off the
coefficient of $u^{2n-1}$ from a {\it polynomial} in $u$ of degree
$j$. Since $j<2n-1$, this coefficient is 0.

If on the other hand we have $i=j$, then (2.7) reduces to
$$\multline (M\cdot U)_{i,j}=
(-1)^i\coef{u^{2n-1}}\frac {1} {1-u}\\
\times\bigg(\frac {1} {\rho(1-\rho+\rho u+xu)}
\bigg(-(1-\rho)(1-\rho u)(1-u)^{i}
\rho^{i-1}\\
+{u(\rho+x)(u^2x+2u-1)}
{(1+xu)^{i-1}} \bigg)\bigg).
\endmultline$$
For the same reasons as before, the term $(1-\rho+\rho u+xu)$ in the
expression in parentheses cancels, whence the expression in
parentheses is in fact a polynomial in $u$ of degree $i+1$. Now we use
the (easily verified) fact that the coefficient of $u^N$ in a formal
power series of the form
$P(u)/(1-q u)$, where $P(u)$ is a polynomial in $u$ of degree $\le N$,
is equal to $q^NP(1/q)$. Thus we obtain
$$(M\cdot U)_{i,j}=\frac {(-1)^i} {\rho}
{(\rho+x)} {(1+x)^{i-1}},
$$
in accordance with (2.4).

\smallskip
{\it Case 4: $j=2n-1$}.
We must evaluate the following sum:
$$\align (M\cdot U)_{i,2n-1}&=\sum _{k=0}
^{2n-1}a_{2n-1-k,i}U_{k,2n-1}\\
&=a_{0,i}+\sum _{k=0} ^{n-1}a_{2n-1-2k,i}U_{2k,2n-1}+
\sum _{k=1} ^{n-1}a_{2n-2k,i}U_{2k-1,2n-1}.
\tag2.8
\endalign$$
Since $U_{0,2n-1}=0$, we may restrict the first sum on the right-hand
side of (2.8) to $1\le k\le n-1$. Consequently, it is equal to
$$\align &-\coef{u^{2n-1}v^i}\sum _{k=1} ^{n-1}
\frac {1} {(\rho+x)(1+x)^{n-2}}\\
&\kern.7cm
\cdot\Bigg(\sum _{s=0} ^{n-k-2}2^{2s+1}
\binom {n-k-1+s}{2s+1} \sum _{t=0} ^{k-s-1}\binom
{n-1}tx^{k-s-t-1}\\
&\kern1.3cm
+\rho \sum _{s=0} ^{n-k-1}2^{2s}
\binom {n-k-1+s}{2s} \sum _{t=0} ^{k-s-1}\binom
{n-1}tx^{k-s-t-1}\Bigg)u^{2k}F(u,v)\\
&=-\coef{u^{2n-1}v^i}
\frac {1} {(\rho+x)(1+x)^{n-2}}\\
&\kern.7cm
\cdot
\Bigg(\sum _{s,k} ^{}2^{2s+1}
\binom {n-k-1}{2s+1} u^{2k+2s}\sum _{t=0} ^{k-1}\binom
{n-1}tx^{k-t-1}\\
&\kern1.3cm
+\rho \sum _{s,k} ^{}2^{2s}
\binom {n-k-1}{2s} u^{2k+2s}\sum _{t=0} ^{k-1}\binom
{n-1}tx^{k-t-1}\Bigg)F(u,v).
\endalign$$
Because of completely analogous arguments, the second sum on the right-hand
side of (2.8) is equal to
$$\align &\coef{u^{2n-1}v^i}\sum _{k=1} ^{n-1}
\frac {1} {(\rho+x)(1+x)^{n-2}}\\
&\kern.7cm
\cdot\Bigg(\sum _{s=0} ^{n-k-1}2^{2s}
\binom {n-k-1+s}{2s} \sum _{t=0} ^{k-s-1}\binom
{n-1}tx^{k-s-t-1}\\
&\kern1.3cm
+\rho \sum _{s=0} ^{n-k-1}2^{2s+1}
\binom {n-k+s}{2s+1} \sum _{t=0} ^{k-s-2}\binom
{n-1}tx^{k-s-t-2}\Bigg)u^{2k-1}F(u,v)\\
&=\coef{u^{2n-1}v^i}
\frac {1} {(\rho+x)(1+x)^{n-2}}\\
&\kern.7cm
\cdot
\Bigg(\sum _{s,k} ^{}2^{2s}
\binom {n-k-1}{2s} u^{2k+2s-1}\sum _{t=0} ^{k-1}\binom
{n-1}tx^{k-t-1}\\
&\kern1.3cm
+\rho \sum _{s,k} ^{}2^{2s+1}
\binom {n-k-1}{2s+1} u^{2k+2s+1}\sum _{t=0} ^{k-1}\binom
{n-1}tx^{k-t-1}\Bigg)F(u,v).
\endalign$$
If this is substituted in (2.8) and the sums are put together, then we obtain
$$\align (M\cdot &U)_{i,2n-1}=a_{0,i}+
\frac {1} {(\rho+x)(1+x)^{n-2}}\coef{u^{2n-1}v^i}(1-\rho u)F(u,v)
\\
&\kern2cm
\times
\sum _{s,k} ^{}(-2)^{s}
\binom {n-k-1}{s} u^{2k+s-1}\sum _{t=0} ^{k-1}\binom
{n-1}tx^{k-t-1}\\
&=a_{0,i}+
\frac {1} {(\rho+x)(1+x)^{n-2}}\coef{u^{2n-1}v^i}(1-\rho u)F(u,v)
\\
&\kern2.5cm
\times
\sum _{k} ^{}(1-2u)^{n-k-1} u^{2k-1}\sum _{t=0} ^{k-1}\binom
{n-1}tx^{k-t-1}\\
&=a_{0,i}+
\frac {1} {(\rho+x)(1+x)^{n-2}}\coef{u^{2n-1}v^i}(1-\rho u)F(u,v)
\\
&\kern2.5cm
\times
\sum _{t=0} ^{n-1}\binom
{n-1}t\frac {1} {u}x^{-t-1}\sum _{k=t+1} ^{n-1}(1-2u)^{n-k-1}
(xu^2)^k\\
&=a_{0,i}+
\frac {1} {(\rho+x)(1+x)^{n-2}}\coef{u^{2n-1}v^i}(1-\rho u)F(u,v)
\\
&\kern2.5cm
\times
\sum _{t=0} ^{n-1}\binom
{n-1}t\frac {1} {u}x^{-t-1}
(1-2u)^{n-1}\frac {\(\frac {xu^2} {1-2u}\)^{t+1}-
\(\frac {xu^2} {1-2u}\)^{n}} {1-\frac {xu^2} {1-2u}}.
\endalign$$
The sum over $t$ is easily evaluated, since it is only geometric
series that have to be summed. After some simplification this leads to
$$\align
(M\cdot &U)_{i,2n-1}=a_{0,i}
+\frac {1} {(\rho+x)(1+x)^{n-2}}\coef{u^{2n-1}v^i}(1-\rho u)\\
&\kern2cm
\times
\frac {u((1-u)^{2n-2}-(u^2(1+x))^{n-1})} {1-2u-xu^2}
F(u,v)\\
&=a_{0,i}
+\frac {1} {(\rho+x)(1+x)^{n-2}}\coef{u^{2n-1}v^i}
\frac {(1-u-v+\rho uv)(u-v)} {(1-u-v-xuv)(1-\rho v)}\\
&\kern2cm
\times
\frac {u((1-u)^{2n-2}-(u^2(1+x))^{n-1})} {1-2u-xu^2}.
\tag2.9
\endalign$$

If we have $i=0$, then it is immediate that the expression (2.9) vanishes.

Hence, from now on we assume $i>0$. We split the expression on the
right-hand side of (2.9) into
$$\align (M&\cdot U)_{i,2n-1}=-\rho^{i-1}
-\frac {1} {(\rho+x)(1+x)^{n-2}}\coef{v^i}
\frac {(1-u-v+\rho uv)(u-v)} {(1-u-v-xuv)(1-\rho v)}\\
&\kern7cm
\times
\frac {(1+x)^{n-1}} {(1-2u-xu^2)}\\
&\kern.9cm
+\frac {1} {(\rho+x)(1+x)^{n-2}}\coef{u^{2n-1}v^i}
\frac {(1-u-v+\rho uv)(u-v)} {(1-u-v-xuv)(1-\rho v)}
\frac {u(1-u)^{2n-2}} {(1-2u-xu^2)}\\
&=-\rho^{i-1}+\rho^{i-1}\frac {1+x} {\rho+x}\\
&\kern.9cm
+\frac {1} {(\rho+x)(1+x)^{n-2}}\coef{u^{2n-1}v^i}
\frac {(1-u-v+\rho uv)(u-v)} {(1-u-v-xuv)(1-\rho v)}
\frac {u(1-u)^{2n-2}} {(1-2u-xu^2)}.\\
\tag2.10
\endalign$$
To the expression in the last line we apply again the partial fraction
expansion (2.6). After having substituted this in (2.10), we may directly read
off the coefficient of $v^i$. Once again, the first term on the
right-hand side of (2.6) does not contribute to it since we have $i>0$
by assumption. We obtain
$$\align (&M\cdot U)_{i,2n-1}\\
&=\rho^{i-1}\frac {1-\rho} {\rho+x}
+\frac {1} {(\rho+x)(1+x)^{n-2}}\coef{u^{2n-1}}
\frac {u(1-u)^{2n-2}} {(1-2u-xu^2)}\frac {1} {(1-\rho+\rho u+xu)}\\
&\kern2cm
\times
\bigg(-{(1-\rho)(1-\rho u)}
\rho^{i-1}
+ {u(\rho+x)(u^2x+2u-1)}
\frac {(1+xu)^{i-1}} {(1-u)^i}\bigg)\\
&=\rho^{i-1}\frac {1-\rho} {\rho+x}
-\frac {1} {(\rho+x)(1+x)^{n-2}}\coef{u^{2n-1}}
\frac {u(1-u)^{2n-2}(1-\rho)(1-\rho u)} {(1-2u-xu^2)(1-\rho+\rho u+xu)}
\rho^{i-1}\\
&\kern1cm
-\frac {1} {(\rho+x)(1+x)^{n-2}}\coef{u^{2n-1}}
  \frac {u^2(1-u)^{2n-2-i}(\rho+x)(1+xu)^{i-1}} {(1-\rho+\rho u+xu)}.
\tag2.11
\endalign$$
Let us for the moment assume that $i<2n-1$.
The coefficient of $u^{2n-1}$ in the last line can be easily
determined by making use of the previously used fact that the
coefficient of $u^N$ in a formal power series of the form
$P(u)/(1-q u)$, where $P(u)$ is a polynomial in $u$ of degree $\le N$,
is equal to $q^NP(1/q)$. Thus we obtain
$$\align (&M\cdot U)_{i,2n-1}\\
&=\rho^{i-1}\frac {1-\rho} {\rho+x}
-\frac {1} {(\rho+x)(1+x)^{n-2}}\coef{u^{2n-1}}
\frac {u(1-u)^{2n-2}(1-\rho)(1-\rho u)} {(1-2u-xu^2)(1-\rho+\rho u+xu)}
\rho^{i-1}\\
&\kern1cm
+\frac {1} {(\rho+x)(1+x)^{n-2}}\frac {(\rho+x)(1+x)^{2n-3}}
{(\rho-1)^{2n-2}}\rho^{i-1}.
\tag2.12
\endalign$$
In order to determine the remaining coefficient of $u^{2n-1}$ on the
right-hand side, we apply once again a partial fraction expansion:
$$\multline
\frac {(1-\rho u)} {(1-2u-xu^2)(1-\rho+\rho u+xu)}=
\frac {\rho+x} {(1+x)(1-\rho+\rho u+xu)}\\
-\frac {1} {2(1+x)}\frac {1}
{(u-\om)}-\frac {1} {2(1+x)}\frac {1}
{(u-\overline\om)},
\endmultline$$
where $\om=-\frac {1} {x}(1-\sqrt{1+x})$ and
$\overline\om=-\frac {1} {x}(1+\sqrt{1+x})$. This is now substituted in
(2.12). By a routine calculation, which uses again the above fact of how
to read off the coefficient of $u^N$ in a formal power series of the form
$P(u)/(1-q u)$, it is seen that (2.12) vanishes for $i<2n-1$.

On the other hand, if we have $i=2n-1$, then the expression in the
last line of (2.11) is equal to
$$-\frac {1} {(\rho+x)(1+x)^{n-2}}\coef{u^{2n-1}}
  \frac {u^2(\rho+x)(1+xu)^{2n-2}} {(1-u)(1-\rho+\rho u+xu)}.$$
In order to read off the coefficient of $u^{2n-1}$ in this expression,
we must apply a partial fraction expansion another time:
$$\frac {1} {(1-u)(1-\rho+\rho u+xu)}=\frac {1} {(1+x)}
\frac {1} {(1-u)}+\frac {1} {(1+x)}\frac {\rho+x}
{(1-\rho+\rho u+xu)}.$$
After having done that, a routine calculation has to be performed that
shows that for $i=2n-1$ the expression (2.11) simplifies to
$-(1+x)^{n-1}$, in accordance with (2.4).

This completes the proof of the theorem.\quad \quad
\qed
\enddemo

Finally we proof the conjectured evaluation of $\det(\tilde A(2n))$ on
page~4 of \cite{\BacRAZ}. As it turns out, this is much simpler than
the previous evaluations.

\proclaim{Theorem 3}Let $(a_{i,j})_{i,j\ge0}$ be the sequence given by the recurrence
$$a_{i,j}=a_{i-1,j}+a_{i,j-1}+x\,a_{i-1,j-1},\quad \quad i,j\ge1
$$
and the initial conditions $a_{i,0}=i$ and $a_{0,i}=-i$,
$i\ge0$. Then
$$\det_{0\le i,j\le 2n-1}(a_{i,j})=(1+x)^{2n(n-1)}.$$
\endproclaim
\demo{Proof}
Again, we compute the generating function
$F(u,v)=\sum _{i,j\ge0} ^{}a_{i,j}u^iv^j$. Here we get
$$F(u,v)=\frac {u-v} {(1-u-v-uvx)(1-u)(1-v)}.$$
This can be written in the following way:
$$\align F(u,v)&=\frac {u} {(1-u)(1-u-v-uvx)}-\frac {v}
{(1-v)(1-u-v-uvx)}\\
&=\frac {u} {(1-u)^2(1-v)\(1-\frac {uv(x+1)} {(1-u)(1-v)}\)}-\frac {v}
{(1-u)(1-v)^2\(1-\frac {uv(x+1)} {(1-u)(1-v)}\)}\\
&=\sum _{\ell=0} ^{\infty}\frac {u^\ell v^\ell
(x+1)^\ell} {(1-u)^{\ell+1}(1-v)^{\ell+1}}\(\frac {u} {1-u}-\frac {v}
{1-v}\).
\endalign$$
Therefore we have
$$a_{i,j}=\sum _{\ell=0} ^{}\(\binom {i}{\ell+1}\binom
{j}\ell-\binom {i}{\ell}\binom
{j}{\ell+1}\)(x+1)^\ell.$$
Now we apply the following row and column operations:
We subtract row $i$ from row $i+1$, $i=2n-2,2n-3,\dots,0$,
and subsequently we subtract
column $j$ from column $j+1$, $j=2n-2,2n-3,\dots,0$. This is repeated
another $2n-2$ times. It is not too difficult to see
(using the above expression for $a_{i,j}$) that, step by
step, the rows and columns are
``emptied'' until finally the determinant
$$\det\pmatrix 0&-1&0&0&\innerhdotsfor2\after\quad &0\\
1&0&-X&0&\innerhdotsfor2\after\quad &0\\
0&X&0&-X^2&\innerhdotsfor2\after\quad &0\\
0&0&X^2&0&\ddots&&0\\
&&\ddots&\ddots&\ddots&\ddots&\vdots\\
0&\innerhdotsfor3\after\quad &X^{2n-3}&0&-X^{2n-2}\\
0&\innerhdotsfor3\after\quad &0&X^{2n-2}&0
\endpmatrix$$
is obtained, where we have written $X$ for $x+1$. The theorem follows
now immediately.\quad \quad \qed
\enddemo

\subhead 3. Some binomial determinants\endsubhead
In this section we prove some evaluations of binomial determinants
which arose from conjectures in \cite{\BacRAZ} and subsequent private
discussion with Roland Bacher.

The first theorem in this section proves a common extension of two
conjectures on page~15 in \cite{\BacRAZ}. It is also
stated (without proof) in \cite{\BacRAA, p.~14}.

\proclaim{Theorem 4}We have
$$\det_{0\le i,j\le n-1}\(\binom {2i+2j+a}i-\binom {2i+2j+a}{i-1}\)
=2^{\binom n2}.\tag3.1$$
\endproclaim
\demo{Proof} We apply LU-factorization. Let us denote the matrix in (3.1)
by $M$.
I claim that
$$M\cdot U=L,$$
where $U=\((-1)^{j-i}\binom ji\)_{0\le i,j\le n-1}$, and where
$L=(L_{i,j})_{0\le i,j\le n-1}$ is a lower triangular matrix with
diagonal entries $1,2,4,8,\dots,2^{n-1}$. It is obvious that the claim
immediately implies the theorem.

Again, we verify the claim by a direct calculation.
Let $0\le i\le j\le n-1$.
The $(i,j)$-entry of $M\cdot U$ is
$$(M\cdot U)_{i,j}=\sum _{k=0} ^{j}
{{\left( -1 \right) }^{j - k}}{j\choose k}
{{\left( 1 + a + 2\,k
  \right) \,\left( a + 2\,i + 2\,k \right) !}\over {i!\,\left( 1
  + a + i + 2\,k \right) !}}.
$$
In standard hypergeometric notation
$${}_r F_s\!\left[\matrix a_1,\dots,a_r\\ b_1,\dots,b_s\endmatrix;
z\right]=\sum _{m=0} ^{\infty}\frac {\po{a_1}{m}\cdots\po{a_r}{m}}
{m!\,\po{b_1}{m}\cdots\po{b_s}{m}} z^m\ ,$$
where the Pochhammer symbol
$(\alpha)_m$ is defined by $(\alpha)_m:=\alpha(\alpha+1)\cdots(\alpha+m-1)$,
$m\ge1$, $(\alpha)_0:=1$, the above binomial sum can be written as
$${{\left( -1 \right) }^j}
{{\left( 1 + a \right) \,
     ({ \textstyle 2 + a + i}) _{i-1} }\over {i!}}
     {} _{4} F _{3} \!\left [ \matrix { {3\over 2} + {a\over 2}, {1\over 2} +
      {a\over 2} + i, 1 + {a\over 2} + i, -j}\\ { {3\over 2} + {a\over 2} +
      {i\over 2}, 1 + {a\over 2} + {i\over 2}, {1\over 2} + {a\over
      2}}\endmatrix ; {\displaystyle 1}\right ].
$$
To this $_4F_3$-series we apply the contiguous relation
$$
{} _{4} F _{3} \!\left [ \matrix { a, b,c,d}\\ {e,f,g}\endmatrix ; {\displaystyle
   z}\right ]  = {} _{4} F _{3} \!\left [ \matrix { a - 1,
b,c,d}\\ {
    e,f,g}\endmatrix ; {\displaystyle z}\right ]  +
   {z    }
   {{bcd }\over
    {efg}}
   {} _{4} F _{3} \!\left [ \matrix { a, b+1,c+1,d+1}\\ {e+1,f+1,g+1}\endmatrix ;
        {\displaystyle z}\right ].
$$
We obtain
$$\multline {{\left( -1 \right) }^j}
{{\left( 1 + a \right) \,
      ({ \textstyle 2 + a + i}) _{i-1} }\over {i!}}
      {} _{3} F _{2} \!\left [ \matrix { {1\over 2} + {a\over 2} + i, 1 +
       {a\over 2} + i,-j}\\ { {3\over 2} + {a\over 2} + {i\over 2}, 1 + {a\over
       2} + {i\over 2}}\endmatrix ; {\displaystyle 1}\right ]\\
-
{{\left( -1 \right) }^j}  {{2j\,
      ({ \textstyle 4 + a + i}) _{i-1} }\over {i!}}
      {} _{3} F _{2} \!\left [ \matrix {  {3\over 2} + {a\over 2} + i, 2
       + {a\over 2} + i, 1 - j}\\ { {5\over 2} + {a\over 2} + {i\over 2}, 2 + {a\over
       2} + {i\over 2}}\endmatrix ; {\displaystyle 1}\right ].
\endmultline$$
Now we apply the transformation formula
(cf\. \cite{\BailAA, Ex.~7, p.~98})
$$
{} _{3} F _{2} \!\left [ \matrix { a, b, -n}\\ { d, e}\endmatrix ;
   {\displaystyle 1}\right ]  =
{{      ({ \textstyle -a - b + d + e}) _{n} }\over {({ \textstyle e}) _{n} }}
{} _{3} F _{2} \!\left [ \matrix { -n, -a + d, -b + d}\\ { d, -a - b + d +
       e}\endmatrix ; {\displaystyle 1}\right ]  ,
\tag3.2$$
where $n$ is a nonnegative integer, to both $_3F_2$-series. The above
expression then becomes
$$\multline
{{\left( -1 \right) }^j}
{{\left( 1 + a \right) \,
({ \textstyle 1 - i}) _{j} \,
      ({ \textstyle 2 + a + i}) _{i-1} }\over
    {i!\,({ \textstyle 1 + {a\over 2} + {i\over 2}}) _{j} }}
      {} _{3} F _{2} \!\left [ \matrix { -j, {1\over 2} - {i\over 2}, 1 -
       {i\over 2}}\\ { {3\over 2} + {a\over 2} + {i\over 2}, 1 - i}\endmatrix
       ; {\displaystyle 1}\right ] \\
-
{{\left( -1 \right) }^j}  {{2j\,
({ \textstyle 1 - i}) _{j-1} \,
      ({ \textstyle 4 + a + i}) _{i-1} }\over
    {i!\,({ \textstyle 2 + {a\over 2} + {i\over 2}}) _{j-1} }}
      {} _{3} F _{2} \!\left [ \matrix { 1 - j, {1\over 2} - {i\over 2}, 1 -
       {i\over 2}}\\ { {5\over 2} + {a\over 2} + {i\over 2}, 1 - i}\endmatrix
       ; {\displaystyle 1}\right ].
\endmultline$$
If $j\ge i$, the first summand is always 0 because of the term
$({ \textstyle 1 - i})
_{j}$. The second summand vanishes as long as $j>i$ because of the term
$({ \textstyle 1 - i}) _{j-1} $. Thus, the matrix $L$ is indeed a
lower triangular matrix. It remains to evaluate the above expression
for $i=j$. As we already mentioned, in that case it reduces to the
second term, which itself reduces to
$$-{{\left( -1 \right) }^i}{{2i\,
({ \textstyle 1 - i}) _{i-1} \,
     ({ \textstyle 4 + a + i}) _{i-1} }\over
   {i!\,({ \textstyle 2 + {a\over 2} + {i\over 2}}) _{i-1} }}
     {} _{2} F _{1} \!\left [ \matrix { {1\over 2} - {i\over 2}, 1 - {i\over
      2}}\\ { {5\over 2} + {a\over 2} + {i\over 2}}\endmatrix ; {\displaystyle
      1}\right ].
$$
This $_2F_1$-series can be evaluated by means of Gau\3'
$_2F_1$-summation (cf\. \cite{\SlatAC, (1.7.6); Appendix (III.3)})
$$
{} _{2} F _{1} \!\left [ \matrix { a, b}\\ { c}\endmatrix ; {\displaystyle
   1}\right ]  = \frac {\Gamma( c)\,\Gamma( c-a-b)} {\Gamma(
c-a)\,\Gamma( c-b )}.
\tag3.3$$
After some further simplification we obtain $2^i$, as desired.\quad \quad \qed
\enddemo

Our next theorem presents a variation of the previous theorem. It is again a
common extension of two conjectures in \cite{\BacRAZ, p.~15}.
It is also stated (without proof) in \cite{\BacRAA, p.~14}.

\proclaim{Theorem 5}We have
$$\det_{0\le i,j\le n-1}\(\binom {2i+2j+a}{i+1}-\binom {2i+2j+a}{i}\)
=2^{\binom n2}\frac {\prod _{i=0} ^{n-1}(a+2i-1)} {n!}.\tag3.4$$
\endproclaim
\demo{Proof} We apply LU-factorization. Let us denote the matrix in
(3.4) by $M$. I claim that
$$M\cdot U=L,$$
where
$$U=\((-1)^{j-i}\binom ji\frac {(a+2j-1)}
{(a+2i-1)}\)_{0\le i,j\le n-1},$$
and where
$L=(L_{i,j})_{0\le i,j\le n-1}$ is a lower triangular matrix with
diagonal entries
$${a-1},\frac {2} {2}(a+1),\frac {4} {3}(a+3),\frac {8} {4}(a+5),\dots,
\frac {2^{n-1}}{n}(a+2n-3).$$
It is obvious that the claim
immediately implies the theorem.

Let us do the required calculations. Let $0\le i\le j\le 2n-1$.
The $(i,j)$-entry of $M\cdot U$ is
$$(M\cdot U)_{i,j}=\sum _{k=0} ^{j}
{{\left( -1 \right) }^{j - k}}{j\choose k}
{{\left( a + 2\,j-1
  \right) \,\left( a + 2\,i + 2\,k \right) !}\over {\left(i+1
  \right) !\,\left( a + i + 2\,k \right) !}}.
$$
We convert this binomial sum into hypergeometric notation. The result is
$${{\left( -1 \right) }^j}{{\left( a + 2\,j -1\right) \,
     (a + 2\,i)! }\over
   {(i+1)!\,(a+i)!}}
     {} _{3} F _{2} \!\left [ \matrix { -j, {1\over 2} + {a\over 2} + i, 1 +
      {a\over 2} + i}\\ { 1 + {a\over 2} + {i\over 2}, {1\over 2} + {a\over 2}
      + {i\over 2}}\endmatrix ; {\displaystyle 1}\right ].
$$
This time it is not even necessary to apply a contiguous relation. We
may directly apply the transformation formula (3.2) and obtain
$${{\left( -1 \right) }^j}
{{\left( a + 2\,j -1\right) \,
(a + 2\,i)! \,
     ({ \textstyle -i}) _{j} }\over
   {(i+1)!\,(a+i)!\,
     ({ \textstyle {1\over 2} + {a\over 2} + {i\over 2}}) _{j} }}
     {} _{3} F _{2} \!\left [ \matrix { -j, -{{i}\over 2}, {1\over 2} -
      {i\over 2}}\\ { 1 + {a\over 2} + {i\over 2}, -i}\endmatrix ;
      {\displaystyle 1}\right ].
$$
Again, this expression is 0 if $j>i$ since it contains the term
$(-i)_j$. If $j=i$ then it reduces to
$${{\left( -1 \right) }^i}
{{\left( a + 2\,i -1\right) \,
(a + 2\,i)! \,({ \textstyle -i}) _{i} }\over
   {(i+1)!\,(a+i)!\,
     ({ \textstyle {1\over 2} + {a\over 2} + {i\over 2}}) _{i} }}
     {} _{2} F _{1} \!\left [ \matrix {- {{i}\over 2}, {1\over 2} - {i\over
      2}}\\ { 1 + {a\over 2} + {i\over 2}}\endmatrix ; {\displaystyle 1}\right
      ].
$$
The $_2F_1$-series can again be evaluated by means of the Gau\3 summation
formula (3.3). After further simplification we obtain
$2^i(a+2i-1)/(i+1)$, as desired.\quad \quad \qed

\enddemo

The final determinant evaluation that we present here proves a common
generalization of determinant evaluations that arose in discussion
with Roland Bacher. It is interesting to note that the limiting case where
$Y\to\infty$, i.e., the evaluation of the determinant $\det_{0\le
i,j\le 2n-1}\((i-j)\,(X+i+j)!  \)$ is covered by a theorem by Mehta
and Wang \cite{\MeWaAA}. In fact, the main theorem of \cite{\MeWaAA}
gives the evaluation of the more general determinant
$\det_{0\le
i,j\le n-1}\((Z+i-j)\,(X+i+j)!  \)$. It seems difficult to find a
common extension of the Mehta and Wang evaluation and the evaluation
below, i.e., an evaluation of
$$\det_{0\le i,j\le n-1}\((Z+i-j)\frac {(X+i+j)!} {(Y+i+j)!} \).$$
Neither the method used by Mehta and Wang nor the ``identification of
factors method" used below seem to help.

\proclaim{Theorem 6}Let $X$ and $Y$ be arbitrary nonnegative integers.
Then
$$\multline
\det_{0\le i,j\le 2n-1}\((i-j)\frac {(X+i+j)!} {(Y+i+j)!} \)
=\prod _{i=0} ^{2n-1}\frac {(X+i)!} {(Y+i+2n-1)!}\\
\times
\prod _{i=0} ^{n-1}(2i+1)!^2 (X+2i+1) (Y+4n-2i-2)
    (X-Y-2i)^{2n-2i-2} (X-Y-2i-1)^{2n-2i-2}.
\endmultline\tag3.5$$
\endproclaim
\remark{Remark} The theorem is formulated only for integral $X$ and
$Y$. But in fact, if we replace the factorials by the appropriate
gamma functions then
the theorem would also make sense and be true for complex $X$ and $Y$.
\endremark
\demo{Proof}
To begin with, we take some factors out of the determinant:
$$\multline
\det_{0\le i,j\le 2n-1}\((i-j)\frac {(X+i+j)!} {(Y+i+j)!} \)=
\prod _{i=0} ^{2n-1}\frac {(X+i)!} {(Y+i+2n-1)!}\\
\times
\det_{0\le i,j\le 2n-1}\big((i-j)\,(X+i+1)_{j}\,(Y+i+j+1)_{2n-j-1} \big).
\endmultline$$
If we compare with (3.5), then we see that it suffices to establish that
$$\multline
\det_{0\le i,j\le 2n-1}\big((i-j)\,(X+i+1)_{j}\,(Y+i+j+1)_{2n-j-1} \big)
\\=
\prod _{a=0} ^{n-1}(2a+1)!^2 (X+2a+1) (Y+4n-2a-2) \\
\cdot
   (X-Y-2a)^{2n-2a-2} (X-Y-2a-1)^{2n-2a-2} .
\endmultline\tag3.6$$

This time we apply a different method. We make use of the
{\it ``identification of factors" method\/} as explained in \cite{\KratBN,
Sec.~2.4}.

Let us denote the determinant in (3.6) by $T_n(X,Y)$.
Here is a brief outline of how we proceed.
The determinant $T_n(X,Y)$
is a polynomial in $X$ and $Y$. In the first step we shall show
that the product on the right-hand side of (3.6) divides the
determinant as a polynomial in $X$ and $Y$. Then we compare the
degrees in $X$ and $Y$ of the determinant and the right-hand side of
(3.6). As it turns out, the degree of the determinant is at most the
degree of the product. Hence, the determinant must be equal to a
constant multiple of the right-hand side of (3.6). Finally, in the
last step, this constant is computed by setting $X=-2n$ on both sides
of (3.6).

\smallskip
{\it Step 1. The product $\prod _{a=0} ^{n-1}(X+2a+1)$ divides
$T_n(X,Y)$}. Let $a$ be fixed with $0\le
a\le n-1$. In order to show that $X+2a+1$ divides the
determinant, it suffices to show that the determinant vanishes for
$X=-2a-1$. The latter would follow immediately if we could show that
the rows of the matrix underlying the determinant are linearly
dependent for $X=-2a-1$. I claim that in fact we have
$$\sum _{i=0} ^{a}(-1)^i\binom ai
\frac {({ \textstyle Y+2n+2a - i  }) _{i}}
{({ \textstyle \frac Y 2+ 2a - i  }) _{i}}
\big(\text {row $(2a-i)$ of
$T_n(-2a-1,Y)$}\big)=0.$$
This claim is easily verified. Indeed, if we restrict it to the $j$-th
column, we see that we must check
$$\sum _{i=0} ^{a}(-1)^i\binom ai
\frac {({ \textstyle Y+2n+2a - i  }) _{i}}
{({ \textstyle \frac Y 2+ 2a - i  }) _{i}}
 {\left( 2\,a - i - j \right)\,
({ \textstyle -i}) _{j} \,({ \textstyle Y + 2a - i + j + 1}) _{2n-j-1}}
=0.
$$
Because of the term $(-i)_j$ this equation is certainly true if
$a<j$. If $a\ge
j$ then, because of the same reason, we may restrict the sum to $j\le i\le
a$. We convert the sum into hypergeometric notation, and obtain
$$2{{({ \textstyle a - j}) _{j+1} \,
     ({ \textstyle 1 + 2\,a + Y}) _{2n-1} }\over
   {({ \textstyle \frac Y 2+2a - j }) _{j} }}
{} _{3} F _{2} \!\left [ \matrix { 1 - 2\,a + 2\,j, -a + j, -2\,a - Y}\\
      { -2\,a + 2\,j, 1 - 2\,a + j - {Y\over 2}}\endmatrix ; {\displaystyle
      1}\right ].$$
Next we apply the contiguous relation
$$
{} _{3} F _{2} \!\left [ \matrix { a, b,c}\\ {d,e}\endmatrix ; {\displaystyle
   z}\right ]  = {} _{3} F _{2} \!\left [ \matrix { a - 1,
b,c}\\ {
   d, e}\endmatrix ; {\displaystyle z}\right ]  +
   {z    }
   {{bc }\over
    {de}}
   {} _{3} F _{2} \!\left [ \matrix { a, b+1,c+1}\\ {d+1,e+1}\endmatrix ;
        {\displaystyle z}\right ].
\tag3.7$$
This yields
$$\multline
2  {{      ({ \textstyle a - j}) _{j+1} \,
      ({ \textstyle 1 + 2\,a + Y}) _{-1 + 2\,n} }\over
    {({ \textstyle 2\,a - j + {Y\over 2}}) _{j} }}
{} _{2} F _{1} \!\left [ \matrix { -2\,a - Y,-a + j}\\ { 1 - 2\,a + j
       - {Y\over 2}}\endmatrix ; {\displaystyle 1}\right ]
\\+
{{
      ({ \textstyle a - j}) _{j+1} \,({ \textstyle 2\,a + Y}) _{2\,n} }\over
      {({ \textstyle -1 + 2\,a - j + {Y\over 2}}) _{j+1} }}
{} _{2} F _{1} \!\left [ \matrix { 1 - 2\,a - Y,1 - a + j}\\ { 2 - 2\,a + j
       - {Y\over 2}}\endmatrix ; {\displaystyle 1}\right ] .
\endmultline$$
Both $_2F_1$-series can be evaluated by means of the Chu-Vandermonde
summation formula (cf\. \cite{\SlatAC, (1.7.7); Appendix (III.4)}
$$
{} _{2} F _{1} \!\left [ \matrix { a, -n}\\ { c}\endmatrix ; {\displaystyle
   1}\right ]  = {{({ \textstyle c-a}) _{n} }\over
    {({ \textstyle c}) _{n} }},
\tag3.8$$
where $n$ is a nonnegative integer. After some further simplification
it is seen that the two terms cancel each other.

\smallskip
{\it Step 2. The product $\prod _{a=0} ^{n-1}(Y+4n-2a-2)$ divides
$T_n(X,Y)$}. We proceed analogously. Let $a$ be fixed with $0\le
a\le n-1$. I claim that
$$\multline
\sum _{i=0} ^{a}(-1)^i\binom ai
\frac {({ \textstyle {X\over 2}+2n - a - i + {1\over 2}}) _{i}}
{({ \textstyle X+2n-a - i }) _{i}}\\
\cdot
\big(\text {row $(2n-a-i-1)$ of
$T_n(X,-4n+2a+2)$}\big)=0.
\endmultline$$
Restricted to the $j$-th column, this is
$$\multline
\sum _{i=0} ^{a}(-1)^i\binom ai
\frac {({ \textstyle {X\over 2}+2n - a - i + {1\over 2}}) _{i}}
{({ \textstyle X+2n-a - i }) _{i}}\\
\cdot
\left( 2n - a - i - j -1 \right)
({\textstyle X+2n-a - i }) _{j}
({ \textstyle a -2n- i + j +2}) _{ 2n-j-1}=0.
\endmultline$$
In order to establish this equation,
we convert the sum again into
hypergeometric notation, and obtain
$$\multline
-  ({ \textstyle 1 + a + j - 2\,n}) _{-j + 2\,n} \,
    ({ \textstyle -a + 2\,n + X}) _{j}  \\
\times
 {} _{3} F _{2} \!\left [ \matrix { 2 + a + j - 2\,n, -1 - a - j +
     2\,n, {1\over 2} + a - 2\,n - {X\over 2}}\\ { 1 + a + j - 2\,n, 1 + a - j
     - 2\,n - X}\endmatrix ; {\displaystyle 1}\right ].
\endmultline$$
Now we apply the contiguous relation (3.7), to get
$$\multline
- ({ \textstyle 1 + a + j - 2\,n}) _{2n-j} \,
   ({ \textstyle -a + 2\,n + X}) _{j} \\
\times
  {} _{2} F _{1} \!\left [ \matrix {  {1\over 2} + a - 2\,n
    - {X\over 2},-1 - a - j + 2\,n}\\ { 1 + a - j - 2\,n - X}\endmatrix ; {\displaystyle
    1}\right ]\\
- \left( \tfrac 1 2 + a - 2\,n - \tfrac X 2 \right) \,
({ \textstyle 1 + a + j - 2\,n}) _{2 n-j} \,
     ({ \textstyle -a + 2\,n + X}) _{j-1}  \\
\times
     {} _{2} F _{1} \!\left [ \matrix {  {3\over 2} + a - 2\,n -
      {X\over 2},-a - j + 2\,n}\\ { 2 + a - j - 2\,n - X}\endmatrix ; {\displaystyle
      1}\right ].
\endmultline$$
Finally, we evaluate the $_2F_1$-series by means of the
Chu-Vandermonde summation formula (3.8).

\smallskip
{\it Step 3. The product $\prod _{a=0} ^{n-1} (X-Y-2a)^{2n-2a-2}$ divides
$T_n(X,Y)$}.
Let $a$ be fixed with $0\le
a\le n-1$. I claim that for $0\le v\le 2n-2a-3$ we have
$$\multline
\sum _{i=0} ^{2a+2}(-1)^i \binom {2a+2}i
\frac {({ \textstyle X +2n+v- i  + 2}) _{i}}
{({ \textstyle X + 2a +v- i + 3}) _{i}}\\
\cdot
\big(\text {row $(2a+v+2-i)$ of
$T_n(X,X-2a)$}\big)=0.
\endmultline$$
Restricted to the $j$-th column, this is
$$\multline
\sum _{i=0} ^{2a+2}(-1)^i \binom {2a+2}i
\frac {({ \textstyle X +2n+v- i  + 2}) _{i}}
{({ \textstyle X + 2a +v- i + 3}) _{i}}\\
\cdot
\left(2a +v- i - j + 2 \right)\, ({
\textstyle X + 2\,a +v - i + 3}) _{j} \,({ \textstyle X+v - i + j +3})
_{ 2n-j-1}=0.
\endmultline$$
In order to establish this equation,
we convert the sum into hypergeometric
notation, and obtain
$$\multline
\left( 2\,a+v - j + 2 \right) \,
 ({ \textstyle X + 2\,a + v + 3}) _{j} \,
  ({ \textstyle X + v + j + 3}) _{2n-j-1} \\
\times
  {} _{3} F _{2} \!\left [ \matrix { -1 - 2a + j - v, -2 - 2a, -2 - j
- v - X}\\ { -2 - 2a + j - v, -2 - 2a - j - v - X}\endmatrix ;
{\displaystyle 1}\right ].
\endmultline$$
Next we apply the contiguous relation (3.7), to get
$$\multline
  \left( 2a+v - j + 2 \right) \,
   ({ \textstyle X + 2\,a + v + 3}) _{j} \,
   ({ \textstyle X + v + j + 3}) _{2n-j-1} \\
\times
   {} _{2} F _{1} \!\left [ \matrix {  -2 - j - v - X,-2 - 2\,a}\\ { -2 -
2\,a - j - v - X}\endmatrix ; {\displaystyle 1}\right ]\\
+
2\left( a +1\right) \,({ \textstyle X + 2a + v + 3}) _{j-1} \,
   ({ \textstyle X + v + j + 2}) _{2n-j}\\
\times
{} _{2} F _{1} \!\left [ \matrix {
-1 - j - v - X,-1 - 2\,a}\\ { -1 - 2a - j - v - X}\endmatrix ; {\displaystyle
1}\right ] .
\endmultline$$
Finally, the $_2F_1$-series are again evaluated by means of the
Chu-Vandermonde summation formula (3.8). The result is
$$\multline
\left( 2 + 2\,a + 3\,j + 2\,a\,j + v + 2\,a\,v + 2\,X + 2\,a\,X
\right)\\
\times
\frac{ ({ \textstyle -2a}) _{2a+1} \,
    ({ \textstyle X + 2a + v + 3}) _{j-1} \,
    ({ \textstyle X + v + j + 3}) _{2n-j-1} }{({ \textstyle  -
2a - j - v - X-1}) _{2a+1} }.
\endmultline$$
Since this expression contains the term $(-2a)_{2a+1}$, it must
vanish.

\smallskip
{\it Step 4. The product $\prod _{a=0} ^{n-1} (X-Y-2a-1)^{2n-2a-2}$ divides
$T_n(X,Y)$}.
Let $a$ be fixed with $0\le
a\le n-1$. I claim that for $0\le v\le 2n-2a-4$ (sic!) we have
$$\multline
\sum _{i=0} ^{2a+3}(-1)^i \binom {2a+3}i
\frac {({ \textstyle X +2n+v- i  + 2}) _{i}}
{({ \textstyle X + 2a +v- i + 4}) _{i}}\\
\cdot
\big(\text {row $(2a+v+3-i)$ of
$T_n(X,X-2a-1)$}\big)=0.
\endmultline$$
Restricted to the $j$-th column, this is
$$\multline
\sum _{i=0} ^{2a+3}(-1)^i \binom {2a+3}i
\frac {({ \textstyle X +2n+v- i  + 2}) _{i}}
{({ \textstyle X + 2a +v- i + 4}) _{i}}\\
\cdot
\left(2a +v- i - j + 3 \right)\, ({
\textstyle X + 2\,a +v - i + 4}) _{j} \,({ \textstyle X+v - i + j +3})
_{ 2n-j-1}=0.
\endmultline$$
In order to establish this equation,
we convert the sum into hypergeometric
notation, and obtain
$$\multline
\left( 2\,a+v - j + 3 \right) \,
 ({ \textstyle X + 2\,a + v + 4}) _{j} \,
  ({ \textstyle X + v + j + 3}) _{2n-j-1} \\
\times
  {} _{3} F _{2} \!\left [ \matrix { -2 - 2a + j - v, -3 - 2a, -2 - j
- v - X}\\ { -3 - 2a + j - v, -3 - 2a - j - v - X}\endmatrix ;
{\displaystyle 1}\right ].
\endmultline$$
Now we apply the contiguous relation (3.7), to get
$$\multline
  \left( 2a+v - j + 3 \right) \,
   ({ \textstyle X + 2\,a + v + 4}) _{j} \,
   ({ \textstyle X + v+j  + 3}) _{2n-j-1}\\
\times
   {} _{2} F _{1} \!\left [ \matrix { -3 - 2a, -2 - j - v - X}\\ { -3 -
2a - j - v - X}\endmatrix ; {\displaystyle 1}\right ]
\\ +
\left(2a+3 \right)  \,({ \textstyle X + 2a + v + 4}) _{j-1} \,
   ({ \textstyle X + v+j + 2}) _{2n-j}  \\
\times
{} _{2} F _{1} \!\left [ \matrix { -2 - 2a,
-1 - j - v - X}\\ { -2 - 2a - j - v - X}\endmatrix ; {\displaystyle
1}\right ].
\endmultline$$
Finally, we evaluate the $_2F_1$-series by means of the
Chu-Vandermonde summation formula (3.8). The result is
$$\multline
\left( 3 + 2\,a + 4\,j + 2\,a\,j + 2\,v + 2\,a\,v + 3\,X +
      2\,a\,X \right)\\
\times
\frac{ ({ \textstyle - 2a-1}) _{2a+2} \,
    ({ \textstyle X + 2a + v + 4}) _{j-1} \,
    ({ \textstyle X + v + j + 3}) _{2n-j-1} }{({ \textstyle -X -
2a - j - v - 2}) _{2a+2} }.
\endmultline$$
Since this expression contains the term $(-2a-1)_{2a+2}$, it must
vanish.

The arguments thus far prove that $\prod _{a=0} ^{n-1} (X-Y-2a-1)^{2n-2a-3}$
divides the determinant $T_n(X,Y)$. It should be observed that the
exponents in this product are by 1 smaller than what we would need.
However, the original determinant in (3.5) is the
determinant of a
skew-symmetric matrix. Hence it is a square (of the corresponding
Pfaffian; see e.g\. \cite{\StemAE, Prop.~2.2}).
This implies that, for fixed $a$, the multiplicity of a factor
$(X-Y-2a-1)$ in $T_n(X,Y)$ must be {\it even}. Phrased differently,
if $(X-Y-2a-1)^e$ divides $T_n(X,Y)$ with $e$ maximal, then $e$ must
be even. Since we already proved that $(X-Y-2a-1)^{2n-2a-3}$ divides
$T_n(X,Y)$, we get for free that in fact $(X-Y-2a-1)^{2n-2a-2}$
divides $T_n(X,Y)$, as required.

\smallskip
{\it Step 5. Comparison of degrees}. It is obvious that the degree in
$X$ of the determinant $T_n(X,Y)$ is at most $\binom {2n}2=2n^2-n$. The
degree in $X$ on the right-hand side of (3.6) is $n+4\binom n2=2n^2-n$,
which is exactly the same number. The same is true for $Y$.
It follows that
$$\multline
T_n(X,Y)=\det_{0\le i,j\le 2n-1}\((i-j)\,(X+i+1)_{j}\,(Y+i+j+1)_{2n-j-1} \)
\\=
C\prod _{a=0} ^{n-1} (X+2a+1) (Y+4n-2a-2) \\
\cdot
   (X-Y-2a-1)^{2n-2a-2} (X-Y-2a)^{2n-2a-2},
\endmultline\tag3.9$$
where $C$ is some constant independent of $X$ and $Y$.

\smallskip
{\it Step 6. The computation of the constant}.
In order to compute the constant, we set $X=-2n$ in (3.9). Because of
$$(-2n+1+i)_j=(-2n+i+1)(-2n+i+2)\cdots(-2n+i+j),$$
this implies that all entries in $T_n(-2n,Y)$ in the
$i$-th row and $j$-th column with $i+j\ge2n$, vanish. This means that the
matrix underlying $T_n(-2n,Y)$ is triangular. Thus its determinant is
easily evaluated. We leave it to the reader to check that the result for $C$
is in accordance with the original claim.\quad \quad \qed
\enddemo


\Refs

\ref\no \BacRAZ\by R.    Bacher
\paper Matrices related to the Pascal triangle
\jour preprint, {\tt math/0109013v1}\vol
\pages \endref

\ref\no \BacRAA\by R.    Bacher \yr 2001
\paper Determinants of matrices related to the Pascal triangle
\jour J. Th\'eorie Nombres Bordeaux\vol 13
\pages \toappear\endref

\ref\no \BailAA\by W. N. Bailey \yr 1935
\book Generalized hypergeometric series
\publ Cambridge University Press
\publaddr Cambridge\endref

\ref\no \KratBN\by C.    Krattenthaler \yr 1999
\paper Advanced determinant calculus
\jour S\'eminaire Lotharingien Combin\.\vol 42 \rm(``The Andrews Festschrift")
\pages Article~B42q, 67~pp\endref

\ref\no \MeWaAA\by M. L. Mehta and R. Wang \yr 2000
\paper Calculation of a certain determinant
\jour Commun\. Math\. Phys\.\vol 214
\pages 227--232\endref

\ref\no \SlatAC\by L. J. Slater \yr 1966
\book Generalized hypergeometric functions
\publ Cambridge University Press
\publaddr Cambridge\endref

\ref\no \StemAE\by J. R. Stembridge \yr 1990
\paper Nonintersecting paths, pfaffians and plane partitions
\jour Adv\. in Math\.\vol 83
\pages 96--131\endref

\endRefs
\enddocument